# Operational Interpretations of the Chernoff Inequality

Roy S. Freedman[1]

**Abstract**

We utilize operational methods to generalize the Chernoff inequality and prove a new result that relates the moment bound to strictly absolute monotonic functions. We show that the Chernoff bound is part of a continuum of probability bounds.

## 1 Introduction

Heaviside [1] developed the operational calculus in the 1890's as a way to solve the differential equations of long-distance telegraphy. Today, most formalisms for the operational method use Laplace and Fourier transforms together with generalized functions (such as the Dirac delta function or the Heaviside step function) as presented by van der Pol [2], Graf [3], and Kanwal [4]. For probability and statistics, the moment generating function and characteristic function corresponds to the Laplace and Fourier transform of the probability density function.

Concentration inequalities provide bounds on the behavior a random variable. Besides being of theoretical interest, these inequalities are useful for bounding the probabilities of random quantities. Our primary objective is to investigate generalizations of the Chernoff inequality [5]. We base our approach on the methods of operational calculus which utilize differential operators that operate on certain classes of functions. A secondary objective is to compare these inequalities with the Chernoff and moment bounds. A consequence of this is to gain an improved understanding of tail probabilities and expected values of functions of random variables.

In Section 2 we provide our notation for operational calculus and probability and review some important results. We discuss our operational interpretation for the Chernoff inequality in Section 3, and show how this generalizes the traditional inequality. Section 4 examines the connection between moment bounds and strictly absolute monotonic functions. Section 5 discusses some preliminary Chernoff-bound results for functions that are not strictly absolutely monotonic.

## 2 Operational Methods and Probability

Heaviside's original approach did not rely on integral transforms. His approach was based on functions and series of differential operators. Mikusinski [6] uses a more algebraic approach that is based on convolution. Our approach here follows the formalism due to Lindell [7] that combines the convolution, differential operator, and transform approaches.

Consider the Taylor series for function $f(z)$. The operational Taylor series for function $f(z)$ is represented by

$$f(z+y) = \sum_{n=0}^{\infty} \frac{y^n}{n!} \cdot \left( \frac{\partial^n}{\partial z^n} \cdot f(z) \right) = \sum_{n=0}^{\infty} \frac{(y \cdot q)^n}{n!} \cdot f(z) = e^{y \cdot q} \cdot f(z) \ . \qquad (2.1)$$

---

[1] Roy S. Freedman is with Inductive Solutions, Inc., New York, NY 10280 and with the Department of Finance and Risk Engineering, New York University Tandon School of Engineering, Brooklyn NY 11201. Email: roy@inductive.com.

Thus, $f(z+y)$ is the result of applying operator $e^{y \cdot q}$ on function $f(z)$. For numerical computation, note that this expression is true for any well-defined $z$. To avoid ambiguity, denote the Taylor series around any point with the $\_|_{z=x}$ notation. For example, the Taylor series around $z = 0$ is

$$f(y) = e^{y \cdot q} \cdot f(z)\big|_{z=0} . \tag{2.2}$$

When we use the operational Taylor series (2.1) in the following convolution integral we obtain:

$$\int_{-\infty}^{\infty} f(z-y) \cdot g(y) \cdot dy = \int_{-\infty}^{\infty} \left(e^{-y \cdot q} \cdot f(z)\right) \cdot g(y) \cdot dy = \left(\int_{-\infty}^{\infty} e^{-y \cdot q} \cdot g(y) \cdot dy\right) \cdot f(z) . \tag{2.3}$$

The convolution corresponds to a differential operator that operates on function $f(z)$. The differential operator is:

$$G(q) = \int_{-\infty}^{\infty} e^{-y \cdot q} \cdot g(y) \cdot dy . \tag{2.4}$$

This is the two-sided Laplace transform. The commutation of convolution implies the commutation of operators:

$$G(q) \cdot f(z) = \int_{-\infty}^{\infty} f(z-y) \cdot g(y) \cdot dy = \int_{-\infty}^{\infty} g(z-w) \cdot f(w) \cdot dw = F(q) \cdot g(z) . \tag{2.5}$$

From the chain rule and (2.5) this implies

$$G(-q) \cdot f(z) = F(q) \cdot g(-z) . \tag{2.6}$$

Other operational results follow when we formally expand the operators themselves in a Taylor series around $q = 0$: for $G_n = \frac{1}{n!} \cdot G^{(n)}(q)\big|_{q=0}$ we have

$$G(q) \cdot f(z) = \left(\sum_{n=0}^{\infty} G_n \cdot q^n\right) \cdot f(z) = \sum_{n=0}^{\infty} G_n \cdot (\partial/\partial z)^n \cdot f(z) = \sum_{n=0}^{\infty} G_n \cdot f^{(n)}(z) . \tag{2.7}$$

For a probabilistic interpretation, suppose random variable $\tilde{z}$ has probability density $p_{\tilde{z}}(z)$. Now use (2.3) in the following convolution integral:

$$\int_{-\infty}^{\infty} f(z-y) \cdot p_{\tilde{z}}(y) \cdot dy = \left(\int_{-\infty}^{\infty} e^{-y \cdot q} \cdot p_{\tilde{z}}(z) \cdot dy\right) \cdot f(z) = P_{\tilde{z}}(q) \cdot f(z) . \tag{2.8}$$

In terms of expectations, this is

$$\mathbf{E}\left[f(z-\tilde{z})\right]=\mathbf{E}\left[\exp(-\tilde{z}\cdot q)\cdot f(z)\right]=\mathbf{E}\left[\exp(-\tilde{z}\cdot q)\right]\cdot f(z)=P_{\tilde{z}}(q)\cdot f(z) \ . \quad (2.9)$$

The expectation shows that the convolution is a differential operator $P_{\tilde{z}}(q)$ that operates on function $f(z)$, where, by (2.7):

$$P_{\tilde{z}}(q)\cdot f(z)=\sum_{n=0}^{\infty}P_n\cdot q^n\cdot f(z)=\sum_{n=0}^{\infty}P_n\cdot f^{(n)}(z) \ . \quad (2.10)$$

Note that the differential operator corresponding to $P_{\tilde{z}}(-q)$ is related to the moment generating function: the Taylor series around $q=0$:

$$P_{\tilde{z}}(q)=\int_{-\infty}^{\infty}e^{-y\cdot q}\cdot p_{\tilde{z}}(y)\cdot dy=\mathbf{E}\left[e^{q\cdot\tilde{z}}\right]=\sum_{n=0}^{\infty}\frac{m_n\cdot(-q)^n}{n!} \ . \quad (2.11)$$

Thus $P_n=(-1)^n\cdot m_n/n!$. One constraint for a probability density is that $P_{\tilde{z}}(0)=m_0=1$.

We introduce generalized functions in a formal sense. Define the impulse or Dirac delta function such that

$$G(q)\cdot\delta(z)=\int_{-\infty}^{\infty}\delta(z-y)\cdot g(y)\cdot dy=\int_{-\infty}^{\infty}\delta(y)\cdot g(z-y)\cdot dy=g(z) \ . \quad (2.12)$$

For $g(z)=1$, (2.12) implies that the impulse is a probability density and is an even function. Define the Heaviside step function: $u(z)=1$ for $z>0$ and $u(z)=0$ for $z<0$. For $g(y)=u(\pm y)$, the above (2.12) implies

$$\left.\begin{array}{l}\int_{-\infty}^{z}\delta(y)\cdot dy=u(z)\\ \int_{z}^{\infty}\delta(y)\cdot dy=u(-z)\end{array}\right\}. \quad (2.13)$$

Since the impulse is even, (2.12) and (2.13) imply that $u(0)=1/2$. Differentiate both sides of (2.13) by applying the operator $q$ to both sides:

$$\left.\begin{array}{l}q\cdot\int_{-\infty}^{z}\delta(y)\cdot dy=\delta(z)=q\cdot u(z)=\frac{\partial}{\partial z}\cdot u(z)\\ q\cdot\int_{z}^{\infty}\delta(y)\cdot dy=-\delta(z)=q\cdot u(-z)=\frac{\partial}{\partial z}\cdot u(-z)\end{array}\right\}. \quad (2.14)$$

Consequently, the inverse operators $U(q)=\frac{1}{q}$ and $U(-q)=\left(\frac{1}{-q}\right)$ corresponds to integration:

$$\left.\begin{aligned} \frac{1}{q} &= \frac{1}{\partial/\partial z} = \int_{-\infty}^{z} \_\cdot dy \\ \left(-\frac{1}{q}\right) &= \frac{1}{-\partial/\partial z} = \int_{z}^{\infty} \_\cdot dy \end{aligned}\right\}. \tag{2.15}$$

This also implies $U(q) \cdot \delta(z) = u(z)$. From (2.5)

$$\left.\begin{aligned} G(q) \cdot u(z) &= U(q) \cdot g(z) = \frac{1}{q} \cdot g(z) = \int_{-\infty}^{z} g(x) \cdot dx \\ G(q) \cdot u(-z) &= U(-q) \cdot g(z) = \left(\frac{1}{-q}\right) \cdot g(z) = \int_{z}^{\infty} g(x) \cdot dx \end{aligned}\right\}. \tag{2.16}$$

The probabilistic interpretation of (2.16) is:

$$\left.\begin{aligned} \Pr[\tilde{z} \leq z] &= \mathbf{E}\bigl[u(z-\tilde{z})\bigr] = P_{\tilde{z}}(q) \cdot u(z) = U(q) \cdot p_{\tilde{z}}(z) \\ \Pr[\tilde{z} \geq z] &= \mathbf{E}\bigl[u(\tilde{z}-z)\bigr] = P_{\tilde{z}}(q) \cdot u(-z) = U(-q) \cdot p_{\tilde{z}}(z) \end{aligned}\right\}. \tag{2.17}$$

In the following sections, it is useful to define, for arbitrary $\varepsilon > 0$, $z^{+} = z + \varepsilon$ and $z^{-} = z - \varepsilon$. Then $u(0^{+}) = 1$ and $u(0^{-}) = 0$.

**Note on Positive Random Variables**

We specify the probability density of a positive random variable $\tilde{z}_{+} > 0$ with the Heaviside unit step function so that

$$P_{\tilde{z}_{+}}(q) \cdot \delta(z) = p_{\tilde{z}_{+}}(z) \cdot u(z) . \tag{2.18}$$

Again, for a probability density $P_{\tilde{z}_{+}}(0) = 1$. It will be useful to consider the positive restriction of an arbitrary density $p_{\tilde{z}}(z)$: denote this restriction by

$$p_{\tilde{z}}^{+}(z) = p_{\tilde{z}}(z) \cdot u(z) . \tag{2.19}$$

The positive restriction is not necessarily a probability density. Here note that

$$P_{\tilde{z}}^{+}(-q) = \int_{0}^{\infty} e^{y \cdot q} \cdot p_{\tilde{z}}(y) \cdot u(y) \cdot dy = \sum_{n=0}^{\infty} \frac{m_{n}^{+} \cdot q^{n}}{n!} . \tag{2.20}$$

All coefficients $m_n^+$ ("positive moments") are positive. We see that the positive restriction $p_{\tilde{z}}^+(z)$ is a probability density only when $P_{\tilde{z}}^+(0) = m_0^+ = \Pr[\tilde{z} > 0] = 1$. Any random variable $\tilde{z}_+ > 0$ conditioned on the event $\tilde{z} > 0$ is induces a density for positive random variable $\tilde{z}_+ > 0$ where

$$p_{\tilde{z}^+}(z) = \frac{p_{\tilde{z}}^+(z)}{\Pr[\tilde{z} > 0]} = \frac{p_{\tilde{z}}^+(z)}{P_{\tilde{z}}^+(0)} = \frac{p_{\tilde{z}}(z) \cdot u(z)}{P_{\tilde{z}}(q) \cdot u(-z)\big|_{z=0^+}} \quad . \tag{2.21}$$

## 3  Operational Chernoff Inequalities

The Chernoff bound generalizes the Markov inequality. In our notation the Markov inequality is:

**Lemma 3.1** (Markov Inequality). For $x > 0$ and positive random variable $\tilde{z}^+$ with density $p_{\tilde{z}^+}(y) = p_{\tilde{z}}(y) \cdot u(y)$:

$$\Pr\left[\tilde{z}^+ \geq x > 0\right] \leq \frac{\mathbf{E}\left[\tilde{z}^+\right]}{x \cdot u(x)} \quad . \tag{3.1}$$

*Proof.*
For fixed $x > 0$, note that $\Pr\left[\tilde{z}^+ \geq x \cdot u(x)\right]$ is a constant. Follow the procedure (2.21) and create a random variable $\tilde{w} \geq x > 0$ conditioned on $\tilde{z}^+ \geq x \cdot u(x)$. The density of $\tilde{w}$ is

$$p_{\tilde{w}}(z) = \frac{p_{\tilde{z}^+}(z)}{\Pr\left[\tilde{z}^+ \geq x \cdot u(x)\right]} \quad . \tag{3.2}$$

Note that $\mathbf{E}[\tilde{w}] \geq x > 0$ so we can write this as $\mathbf{E}[\tilde{w}] \geq x \cdot u(x)$. Now use (3.2):

$$\mathbf{E}[\tilde{w}] = \mathbf{E}\left[\frac{\tilde{z}^+}{\Pr\left[\tilde{z}^+ \geq x\right]}\right] = \frac{\mathbf{E}\left[\tilde{z}^+\right]}{\Pr\left[\tilde{z}^+ \geq x\right]} \geq x \cdot u(x).$$

Untangling the probability in the fraction yields inequality (3.1).   •

Note that for any positive function $h(z) \geq 0$, we can similarly create a positive random variable $\tilde{z}^+ = h(\tilde{z})$ so that:

$$\Pr\left[h(\tilde{z}) \geq x > 0\right] \leq \frac{\mathbf{E}[h(\tilde{z})]}{x \cdot u(x)} \quad . \tag{3.3}$$

**Theorem 3.1.** (Operational Chernoff Inequality). Given an arbitrary random variable $\tilde{z}$ and $f(z) \geq 0$ a non-decreasing positive function. Then for $x > 0$ and $f(x+z) > 0$

$$\Pr[\tilde{z} \geq x] \leq \frac{\mathbf{E}[f(\tilde{z}+z)]}{f(x+z)} = \frac{P_{\tilde{z}}(-q) \cdot f(z)}{e^{q \cdot x} \cdot f(z)} = \frac{F(q) \cdot p_{\tilde{z}}(-z)}{e^{q \cdot x} \cdot f(z)} \quad . \tag{3.4}$$

Note that this inequality is a function of real variable $z$: it corresponds to a set of inequalities depending on a point $z$ and a function $f(z)$.

*Proof.* Since $f(z)$ is positive and non-decreasing,

$$\Pr[\tilde{z} \geq x] = \Pr[\tilde{z}+z \geq x+z] = \Pr[f(z+\tilde{z}) \geq f(z+x)].$$

Next, by the operational Taylor series and definition of convolution from (2.6):

$$\mathbf{E}[f(z+\tilde{z})] = \mathbf{E}[e^{q \cdot \tilde{z}} \cdot f(z)] = \mathbf{E}[e^{q \cdot \tilde{z}}] \cdot f(z) = P_{\tilde{z}}(-q) \cdot f(z) = F(q) \cdot p_{\tilde{z}}(-z).$$

Thus

$$\Pr[\tilde{z} \geq x] \leq \frac{\mathbf{E}[f(z+\tilde{z})]}{f(z+x)} = \frac{P_{\tilde{z}}(-q) \cdot f(z)}{e^{q \cdot x} \cdot f(z)} \quad . \quad \square$$

**Definition 3.1.** The Operational Chernoff Bound is the optimal point $z^*$ and optimal function $f^*(z)$ such that

$$\Pr[\tilde{z} \geq x] \leq \frac{P_{\tilde{z}}(-q) \cdot f(z)\big|_{\substack{f=f^* \\ z=z^*}}}{e^{q \cdot x} \cdot f(z)\big|_{\substack{f=f^* \\ z=z^*}}} = \min_{z,f} \frac{P_{\tilde{z}}(-q) \cdot f(z)}{e^{q \cdot x} \cdot f(z)} \quad . \tag{3.5}$$

Usually the function or point or both is unspecified. We can find an Operational Bound via operational series or with convolutions. Here are some examples that illustrate the inequality.

**Example 3.1. Traditional Chernoff Bound**. Consider the class of exponential functions with positive exponent $f(z) = e^{\alpha \cdot z}$ with $\alpha > 0$. Note that $f(z)$ is a positive non-decreasing function. Thus by operational series (2.11):

$$\begin{aligned} P_{\tilde{z}}(-q) \cdot e^{\alpha \cdot z} &= \sum_{n=0}^{\infty} \frac{m_n \cdot q^n}{n!} \cdot e^{\alpha \cdot z} = \sum_{n=0}^{\infty} \frac{m_n \cdot \alpha^n}{n!} \cdot e^{\alpha \cdot z} = e^{\alpha \cdot z} \cdot \sum_{n=0}^{\infty} \frac{m_n \cdot \alpha^n}{n!} \\ &= e^{\alpha \cdot z} \cdot P_{\tilde{z}}(-\alpha) \quad . \end{aligned} \tag{3.6}$$

Theorem 3.1 implies

$$\Pr[\tilde{z} \geq x] \leq \frac{e^{\alpha \cdot z} \cdot P_{\tilde{z}}(-\alpha)}{e^{\alpha \cdot x} \cdot e^{\alpha \cdot z}} = \frac{P_{\tilde{z}}(-\alpha)}{e^{\alpha \cdot x}} \quad . \tag{3.7}$$

Since terms involve $z$ factor out, the inequality is true for all $z$. For the Chernoff Bound, we find the best $\alpha^*$ (this corresponds to the best $f^*$) that minimizes

$$\Pr[\tilde{z} \geq x] \leq \min_{\alpha} \frac{P_{\tilde{z}}(-\alpha)}{e^{\alpha \cdot x}} \quad . \tag{3.8}$$

For example, for the standard normal density, $P_{\tilde{z}}(-q) = \exp(q^2/2)$, the $\alpha$ that minimizes $\exp(-\alpha \cdot x + \alpha^2/2)$ is $\alpha^* = -x$. The Chernoff Bound is then $\exp(-x^2/2)$. ☐

**Example 3.2. The Heaviside-Chernoff Bound.**
For $f(z) = u(z)$

$$\Pr[\tilde{z} \geq x] \leq \min_{z} \frac{P_{\tilde{z}}(-q) \cdot u(z)}{u(z+x)} \quad .$$

The numerator is the convolution

$$P_{\tilde{z}}(-q) \cdot u(z) = \int_{-\infty}^{\infty} p_{\tilde{z}}(-y) \cdot u(z-y) \cdot dy = \int_{-\infty}^{\infty} p_{\tilde{z}}(y) \cdot u(z+y) \cdot dy = \int_{-z}^{\infty} p_{\tilde{z}}(y) \cdot u(z+y) \cdot dy \, .$$

The numerator decreases when $-z$ increases. However, the Heaviside function in the denominator implies that $u(z+x) = 1$ when $z + x = \varepsilon > 0$. The largest $-z$ can be is $-z \leq x^-$. The denominator is zero when $-z \geq x^+$. Thus the quotient is minimized when $z = -x^+$. The inequality yields

$$\Pr[\tilde{z} \geq x] \leq \int_{x^-}^{\infty} p_{\tilde{z}}(y) \cdot dy \quad .$$

For continuous densities, the expression on the right side of the inequality is the exact expression for the probability that $\tilde{z} \geq x$. We call this the Heaviside-Chernoff bound. Recall the probability from (2.17)

$$\Pr[\tilde{z} \geq x] = \mathbf{E}[u(\tilde{z} - x)] = P_{\tilde{z}}(q) \cdot u(-z)\big|_{z=x} = \frac{P_{\tilde{z}}(-q) \cdot u(z)}{u(z+x)}\bigg|_{z=-x^+} \quad . \; ☐ \tag{3.9}$$

Do functions that approximate the Heaviside step function provide the tightest probabilistic bounds? This will be discussed in Section 4.

**Example 3.3. The Moment Bound.** Consider the non-decreasing function $f(z) = z^\alpha \cdot u(z)$. Note that convolution (2.8) implies

$$P_{\tilde{z}}(-q) \cdot \left(z^\alpha \cdot u(z)\right)\Big|_{z=0} = \int_{-\infty}^{\infty} \left((y+z)^\alpha \cdot u(y+z)\right) \cdot p_{\tilde{z}}(y) \cdot dy \Big|_{z=0}$$
$$= \int_{0}^{\infty} y^\alpha \cdot p_{\tilde{z}}(y) \cdot dy = m_\alpha^+ \quad . \qquad (3.10)$$

Evaluating the Operational Chernoff Bound (3.5) at $z = 0$ specifies the moment bound:

$$\Pr[\tilde{z} \geq x] \leq \min_\alpha \frac{P_{\tilde{z}}(-q) \cdot \left(z^\alpha \cdot u(z)\right)\Big|_{z=0}}{(z+x)^\alpha \cdot u(z+x)\Big|_{z=0}} = \min_\alpha \frac{m_\alpha^+}{x^\alpha \cdot u(x)} \quad . \quad \square \qquad (3.11)$$

Are there other non-decreasing functions $f(z)$ that yield tighter bounds than the moment bound? In the next section we show that, for the class of strictly absolutely monotonic functions, the moment bound is the tightest.

## 4    Strictly Absolutely Monotonic Functions and Moment Bounds

Recall from Widder [8] and Feller [9], a function is strictly absolutely monotonic at $z$ if it and all its derivatives are strictly positive at $z$. For example, the exponential function $\exp(\alpha \cdot z)$ is strictly absolutely monotonic for all $z$ when $\alpha > 0$. Other examples of strictly absolutely monotone functions are the classical Mittag-Leffler functions.

Look at positive functions $f(z) \geq 0$ non-decreasing on $(-\infty, \infty)$ and its positive restriction $f^+(z) := f(z) \cdot u(z)$, where in this case, $f(z)$ is positive and non-decreasing on the semi-infinite interval $(0, \infty)$. Similarly, consider density $p_{\tilde{z}}(z) = P_{\tilde{z}}(q) \cdot \delta(z)$ defined on $(-\infty, \infty)$ and its positive restriction defined on $(0, \infty)$ with $p_{\tilde{z}}^+(y) \cdot u(y) := P_{\tilde{z}}^+(-q) \cdot \delta(z)$. Consider the following four convolutions:

$$\left.\begin{array}{ll} (a) & P_{\tilde{z}}(-q) \cdot f(z) = \int_{-\infty}^{\infty} p_{\tilde{z}}(y) \cdot f(z+y) \cdot dy \\[6pt] (b) & P_{\tilde{z}}^+(-q) \cdot f(z) = \int_{0}^{\infty} p_{\tilde{z}}(y) \cdot f(z+y) \cdot dy \\[6pt] (c) & P_{\tilde{z}}(-q) \cdot f^+(z) = \int_{-z}^{\infty} p_{\tilde{z}}(y) \cdot f(z+y) \cdot dy \\[6pt] (d) & P_{\tilde{z}}^+(-q) \cdot f^+(z) = \int_{\max(0,-z)}^{\infty} p_{\tilde{z}}(y) \cdot f(z+y) \cdot dy \end{array}\right\} . \qquad (4.1)$$

**Lemma 4.1.** Given above convolutions and for $f(z) \geq 0$ non-decreasing. For all $-z \geq 0$:
$(c) = (d) \leq (b) \leq (a)$.

*Proof.* The area under the infinite integral $(-\infty, \infty)$ is greater than the area under the semi-infinite integral $(-z, \infty)$ for positive non-decreasing functions. □

Now recall the Cauchy inequality for series — sometimes called "Cauchy's Third Inequality" — see Mitrinović and Vasić [10], Steele [11]:

**Lemma 4.2.** For positive real numbers $a_n$, $b_n$ and sums $S_a = \sum_{n=0}^{N} a_n$ and $S_b = \sum_{n=0}^{N} b_n$

$$r = \min_n \frac{a_n}{b_n} \leq \frac{\sum_{n=0}^{N} a_n}{\sum_{n=0}^{N} b_n} \leq \max_n \frac{a_n}{b_n} = R \ . \quad (4.2)$$

*Proof.* Note that $r \cdot b_n \leq a_n \leq R \cdot b_n$. Summing both sides over $n$ yields

$$r \cdot \sum_{n=0}^{N} b_n \leq \sum_{n=0}^{N} a_n \leq R \cdot \sum_{n=0}^{N} b_n$$

Dividing both sides of the inequalities by $S_b$ gives the result. Note that for $N \to \infty$, the proof is the same as long as $S_a$ and $S_b$ are convergent. □

**Lemma 4.3.** If $f(z)$ is strictly absolutely monotonic then

$$\min_n \frac{m_n^+}{x^n} \leq \min_{z,f} \frac{P_{\tilde{z}}^+(-q) \cdot f(z)}{f(x+z)} \ . \quad (4.3)$$

*Proof.* Express in series:

$$\frac{P_{\tilde{z}}^+(-q) \cdot f(z)}{f(x+z)} = \frac{\sum_{n=0}^{\infty} \frac{m_n^+}{n!} \cdot f^{(n)}(z)}{\sum_{n=0}^{\infty} \frac{x^n}{n!} \cdot f^{(n)}(z)}.$$

Since $f(z)$ is strictly absolutely monotonic, for all $n$, $f^{(n)}(z) > 0$. The conditions for Cauchy's Third Inequality (Lemma 4.2) are satisfied with $a_n = m_n^+ \cdot f^{(n)}(z)/n!$ and $b_n = x^n \cdot f^{(n)}(z)/n!$. Since this lower bound is true for all $z$ and for all strictly absolutely monotone functions, it is true for the smallest $(z^*, f^*(z))$ that minimizes the quotient. □

**Theorem 4.1**. For $f(z)$ any strictly absolutely monotonic function, $\tilde{z}$ an arbitrary random variable and $-z \geq 0$:

$$\Pr[\tilde{z} \geq x] \leq \min_{\alpha} \frac{m_{\alpha}^{+}}{x^{\alpha}} \leq \min_{z,f} \frac{P_{\tilde{z}}(-q) \cdot f(z)}{f(x+z)} \ . \tag{4.4}$$

*Proof.* For arbitrary random variables, as long as $-z \geq 0$, Lemma 4.1 shows

$$P_{\tilde{z}}^{+}(-q) \cdot f(z) \leq P_{\tilde{z}}(-q) \cdot f^{+}(z) \leq P_{\tilde{z}}(-q) \cdot f(z) \ .$$

Dividing by $f(x+z)$

$$\frac{P_{\tilde{z}}^{+}(-q) \cdot f(z)}{f(x+z)} \leq \frac{P_{\tilde{z}}(-q) \cdot f^{+}(z)}{f(x+z)} \leq \frac{P_{\tilde{z}}(-q) \cdot f(z)}{f(x+z)} \ .$$

Theorem 3.1, (3.5), and Lemma 4.3 imply

$$\Pr[\tilde{z} \geq x] \leq \min_{n} \frac{m_{n}^{+}}{x^{n}} \leq \min_{z,f} \frac{P_{\tilde{z}}(-q) \cdot f(z)}{f(x+z)} \ .$$

Instead of minimizing over all positive integers $n \geq 0$, now minimize over all positive reals $\alpha \geq 0$. Since the set of positive real numbers contain the positive integers, we have

$$\Pr[\tilde{z} \geq x] \leq \min_{\alpha} \frac{m_{\alpha}^{+}}{x^{\alpha}} \leq \min_{n} \frac{m_{n}^{+}}{x^{n}} \ . \qquad \square$$

In some sense. the moment bound is the best we can do with strictly absolutely monotonic increasing functions. Theorem 4.1 shows that the moment bound is the tighter than any operational Chernoff bound applied to any strictly absolutely monotonic function. For the strictly absolutely monotonic exponential function $f(z) = \exp(\alpha \cdot z)$, Theorem 4.1 and (3.6) imply

$$\Pr[\tilde{z} \geq x] \leq \min_{\alpha} \frac{m_{\alpha}^{+}}{x^{\alpha}} \leq \min_{\alpha} \frac{P_{\tilde{z}}(-\alpha)}{e^{\alpha \cdot x}} \ .$$

The moment bound is tighter than the traditional Chernoff bound. This was first proven by Philips and Nelson [12]. Theorem 4.1 generalizes their result to all $f(z)$ that are strictly absolutely monotonic functions.

## 5   Discussion: Functions Not Strictly Absolutely Monotonic

One area of future work on the Operational Chernoff Inequality is to find easily computable asymptotic expansions for functions that are not strictly absolutely monotonic. One approach is to consider functions $u_{\alpha}(z)$ and $U_{\alpha}(q)$, with $u_{\alpha}(z) = U_{\alpha}(q) \cdot \delta(z)$ that approximate the Heaviside step function. Operationally this implies that as $\alpha \to 0$ then

$$u_\alpha(z) \to u(z) \text{ and } U_\alpha(q) \to \frac{1}{q} \ .$$

We can use these results and operational series expansions (2.7) to derive some interesting operational identities. From (2.17),

$$\Pr[\tilde{z} \geq z] = \lim_{\alpha \to 0} \mathbf{E}\left[u_\alpha(\tilde{z} - z)\right] = \lim_{\alpha \to 0} P_{\tilde{z}}(q) \cdot u_\alpha(-z) = \lim_{\alpha \to 0} U_\alpha(-q) \cdot p_{\tilde{z}}(z) \ . \quad (5.1)$$

We can also use the Operational Chernoff Bound directly on the Heaviside function approximation. From (3.5),

$$\Pr[\tilde{z} \geq x] \leq \min_\alpha \frac{P_{\tilde{z}}(-q) \cdot u_\alpha(z)\big|_{z=z^*}}{e^{q \cdot x} \cdot u_\alpha(z)\big|_{z=z^*}} = \min_\alpha \frac{U_\alpha(q) \cdot p_{\tilde{z}}(-z)\big|_{z=z^*}}{e^{q \cdot x} \cdot u_\alpha(z)\big|_{z=z^*}} \ . \quad (5.2)$$

Preliminary review shows that most of these operational identities are not practical: the resultant series oscillate or converge slowly if at all.

**Example 5.1.** $f(z, \alpha) = e^{\alpha \cdot z} \cdot u(z)$ . This function is positive, increasing, but not strictly absolutely monotonic since the function and its derivatives vanish for $z < 0$. It is not continuous at zero. Note that

$$U_\alpha(z) = \frac{1}{q - \alpha} = \begin{cases} \dfrac{1}{q} \cdot \dfrac{1}{1 - \alpha/q} = \sum_{n=0}^\infty \dfrac{\alpha^n}{q^{n+1}} \\ \dfrac{1}{-\alpha} \cdot \dfrac{1}{1 - q/\alpha} = \dfrac{1}{-\alpha} \cdot \sum_{n=0}^\infty \dfrac{q^n}{\alpha^n} \end{cases} \ .$$

The first $P_{\tilde{z}}(q) \cdot u_\alpha(-z)$ is a series expressed in terms of moments and derivatives of impulse functions (from the higher derivatives of the Heaviside step function). The second series $U_\alpha(-q) \cdot p_{\tilde{z}}(z)$ results either in successive integrations of the density (which does not help) or in a series containing derivatives of the given density function – a sort of Gram-Charlier series [13]. For the Gaussian density, these types of series are known to oscillate wildly.

**Example 5.2.** $u_\alpha(z) = z^\alpha \cdot u(z)$. This increasing function is not strictly absolutely monotonic since the function and its derivatives vanish for $z < 0$ and for $z > 0$, $f^{(k)}(z) = 0$ for $k > n$. It is continuous at zero since $f(0^+) = 0$. The Operational Bound is not defined for $-z = x^+$ but is defined for other values $-z > x^+$. Note that

$$U_\alpha(z) = \frac{\Gamma(\alpha + 1)}{q^{\alpha+1}} \ .$$

The first series $P_{\tilde{z}}(q) \cdot u_\alpha(-z)$ results in a series in terms of moments and derivatives of impulse functions (from the higher derivatives of the Heaviside step function). The second series $U_\alpha(-q) \cdot p_{\tilde{z}}(z)$ results in successive integrations of the density (which does not help). Note that using the convolution integral (3.10) explicitly results in the moment bound.

**Example 5.3**. $u_\alpha(z) = 1/(1 + e^{-z/\alpha})$. This is the logistic function. This function is positive, increasing, but not strictly absolutely monotonic since the derivatives alternate in sign. Operationally, the logistic function approximates the Heaviside step function since

$$U_\alpha(q) = \frac{\pi \cdot \alpha}{\sin(\pi \cdot \alpha \cdot q)} = \frac{1}{q} + \frac{\alpha^2 \cdot \pi^2 \cdot q}{6} + \frac{7 \cdot \alpha^4 \cdot \pi^4 \cdot q^3}{360} + \frac{31 \cdot \alpha^6 \cdot \pi^6 \cdot q^5}{15120} + ... \quad .$$

The first series $P_{\tilde{z}}(q) \cdot u_\alpha(-z)$ results in an interesting series in terms of moments and derivatives of the logistic function. We can obtain Chernoff-type bounds for certain $z = z^*$ and $\alpha = \alpha^*$. The second series $U_\alpha(-q) \cdot p_{\tilde{z}}(z)$ results in an integration of the density (which does not help) together with a Gram-Charlier type series.